\documentstyle[11pt,amssymb,amsfonts]{article}

\begin{document}
\newcommand{\p}{\parallel }
\makeatletter \makeatother
\newtheorem{th}{Theorem}[section]
\newtheorem{lem}{Lemma}[section]
\newtheorem{de}{Definition}[section]
\newtheorem{rem}{Remark}[section]
\newtheorem{cor}{Corollary}[section]
\renewcommand{\theequation}{\thesection.\arabic {equation}}

\title{{\bf A Kastler-Kalau-Walze type theorem and the spectral action for perturbations of Dirac operators on manifolds with boundary}
}
\author{Yong Wang \\}

\date{}
\maketitle

\begin{abstract} In this paper, we prove a Kastler-Kalau-Walze type theorem for perturbations of Dirac operators on compact manifolds with or without boundary. As a corollary, we give
two kinds of operator-theoretic explanations of the gravitational
action on boundary.
 We also compute the spectral action for Dirac operators with two-form perturbations on $4$-dimensional compact manifolds.\\

\noindent{\bf Keywords:}\quad
Perturbations of Dirac operators; noncommutative residue; gravitational action; spectral action; Seeley-deWitt coefficients\\
\end{abstract}

\section{Introduction}
    \quad The noncommutative residue found in [Gu] and [Wo] plays a
prominent role in noncommutative geometry. In [Co1], Connes used the
noncommutative residue to derive a conformal 4-dimensional Polyakov
action analogy. In [Co2], Connes proved that the noncommutative
residue on a compact manifold $M$ coincided with the Dixmier's trace
on pseudodifferential operators of order $-{\rm {dim}}M$. Several
years ago, Connes made a challenging observation that the
noncommutative residue of the square of the inverse of the Dirac
operator was proportional to the Einstein-Hilbert action, which is
called the Kastler-Kalau-Walze theorem now. In [Ka], Kastler gave a
brute-force proof of this theorem. In [KW], Kalau and Walze proved
this theorem
by the normal coordinates way simultaneously. In [Ac], Ackermann gave a note on a new proof of this theorem by the heat kernel expansion way.
The Kastler-Kalau-Walze theorem had been generalized to some cases, for example, Dirac operators with torsion [AT], CR manifolds [Po], ${\mathbb{R}}^n$[BC1]
(see also [BC2], [Ni]).\\
\indent On the other hand, Fedosov et al. defined a noncommutative
residue on Boutet de Monvel's algebra and proved that it was the
unique continuous trace in [FGLS]. In [Sc], Schrohe gave the relation
between the Dixmier trace and the noncommutative residue for
manifolds with boundary. In [Wa1],[Wa2], we gave an
operator theoretic explaination of the gravitational action for
manifolds with boundary and proved a Kastler-Kalau-Walze type theorem
for Dirac operators and signature operators on manifolds with
boundary.\\
   \indent Perturbations of Dirac operators were investigated by
   several authors.  In [SZ], Sitarz and Zajac investigated the spectral action for scalar perturbations of Dirac operators.
 In [IL,p.305], Iochum and Levy computed the heat kernel coefficients for Dirac operators with one-form perturbations. In [HPS], Hanisch, Pf\"{a}ffle and
    Stephan derived a formula for the gravitational part of the
    spectral action for Dirac operators on $4$-dimensional spin
    manifolds with totally anti-symmetric torsion and this is a
    perturbation with three-forms of Dirac operators. On the other
    hand, in [CM], Connes and Moscovici considered the conformal
    perturbations of Dirac operators. Investigating the perturbations of Dirac operators has some significance (see [IL],[HPS],[CC]).
 Motivated by
    [IL], [SZ] and [HPS], we study the Dirac operators with general form perturbations. We
    prove a Kastler-Kalau-Walze type theorem for general form perturbations and the conformal perturbations of Dirac operators for compact manifolds with or without
     boundary.
 We also compute the spectral action for Dirac operators with two-form perturbations on $4$-dimensional compact manifolds and give detailed computations of
 spectral action for scalar perturbations of Dirac operators in [SZ].\\
      \indent This paper is organized as follows. In Section 2, we
      prove the Lichnerowicz formula for perturbations of Dirac
      operators and prove a Kastler-Kalau-Walze type theorem for perturbations of Dirac operators on $4$-dimensional compact manifolds with or without
      boundary. In Section 3, we prove a Kastler-Kalau-Walze type theorem for conformal perturbations of Dirac operators on compact manifolds with or without
      boundary.
      In Section 4, We compute the spectral action for Dirac operators with scalar and two-form perturbations on $4$-dimensional compact
      manifolds.\\

\section{A Kastler-Kalau-Walze type theorem for perturbations of Dirac operators }

\noindent{\bf 2.1 A Kastler-Kalau-Walze type theorem for
perturbations of Dirac operators on manifolds without boundary}\\

 \indent Let $M$ be a smooth compact Riemannian $n$-dimensional
 manifold without boundary and $V$ be a vector bundle on $M$. Recall that a differential operator $P$
 is of Laplace type if it has locally the form
 $$P=-(g^{ij}\partial_i\partial_j+A^i\partial_i+B),\eqno(2.1)$$
 where $\partial_i$ is a natural local frame on $TM$ and $g_{i,j}=g(\partial_i,\partial_j)$ and
 $(g^{ij})_{1\leq i,j\leq n}$ is the inverse matrix associated
 to the metric matrix $(g_{i,j})_{1\leq i,j\leq n}$ on $M$, and $A^i$ and $B$ are smooth sections of
 ${\rm End}(V)$ on $M$ (endomorphism). If $P$ is a Laplace type operator
 of the form (2.1), then (see [Gi]) there is a unique connection
 $\nabla$ on $V$ and an unique endomorphism $E$ such that
 $$P=-[g^{ij}(\nabla_{\partial_i}\nabla_{\partial_j}-\nabla_{\nabla^L_{\partial_i}{\partial_j}})+E],\eqno(2.2)$$
 where $\nabla^L$ denotes the Levi-civita connection on $M$.
 Moreover (with local frames of $T^*M$ and $V$), $\nabla_{\partial_i}=
 \partial_i+\omega_i$ and $E$ are related to $g^{ij},~A^i$ and $B$
 through
 $$\omega_i=\frac{1}{2}g_{ij}(A^j+g^{kl}\Gamma_{kl}^j{\rm
 Id}),\eqno(2.3)$$
 $$E=B-g^{ij}(\partial_i(\omega_j)+\omega_i\omega_j-\omega_k\Gamma^k_{ij}),\eqno(2.4)$$
where $\Gamma_{ij}^k$ are the Christoffel coefficients of
$\nabla^L.$\\
\indent Now we let $M$ be a $n$-dimensional oriented spin manifold
with Riemannian metric $g$. We recall that the Dirac operator $D$ is
locally given as follows in terms of orthonormal frames $e_i,~1\leq
i\leq n$ and natural frames $\partial_i$ of $TM$: one has
$$D=\sum_{i,j}g^{ij}c(\partial_i)\nabla^S_{\partial_j}=\sum_{i}c(e_i)\nabla^S_{e_i},\eqno(2.5)$$
where $c(e_i)$ denotes the Clifford action which satisfies the relation\\
 $$c(e_i)c(e_j)+c(e_j)c(e_i)=-2\delta_i^j,$$ and
$$\nabla^S_{\partial_i}=\partial_i+\sigma_i,~~\sigma_i=\frac{1}{4}\sum_{j,k}\left<\nabla^L_{\partial_i}e_j,e_k\right>c(e_j)c(e_k).\eqno(2.6)$$
Let
$$\partial^j=g^{ij}\partial_i,~~\sigma^i=g^{ij}\sigma_j,~~\Gamma^k=g^{ij}\Gamma_{ij}^k.\eqno(2.7)$$
By (6a) in [Ka], we have
$$D^2=-g^{ij}\partial_i\partial_j-2\sigma^j\partial_j+\Gamma^k\partial_k-g^{ij}[\partial_i(\sigma_j)+\sigma_i\sigma_j-\Gamma_{ij}^k\sigma_k]+\frac{1}{4}s,
\eqno(2.8)$$ where $s$ is the scalar curvature. Let $\Psi$ be a
smooth differential form on $M$ and we also denote the associated
Clifford action by $\Psi$. We will compute $D_\Psi^2:=(D+\Psi)^2$.
We note that
$$(D+\Psi)^2=D^2+D\Psi+\Psi D+\Psi^2,\eqno(2.9)$$
$$D\Psi+\Psi D=\sum_{ij}g^{ij}\left(c(\partial_i)\Psi+\Psi
c(\partial_i)\right)\partial_j+\sum_{ij}g^{ij}\left(c(\partial_i)\partial_j(\Psi)+c(\partial_i)\sigma_j\Psi+\Psi
c(\partial_i)\sigma_j\right).\eqno(2.10)$$ By (2.8)-(2.10), we have
$$D_\Psi^2=-g^{ij}\partial_i\partial_j+\left(-2\sigma^j+\Gamma^j+c(\partial^j)\Psi+\Psi
c(\partial^j)\right)\partial_j$$
$$+g^{ij}[-\partial_i(\sigma_j)-\sigma_i\sigma_j+\Gamma_{ij}^k\sigma_k
+c(\partial_i)\partial_j(\Psi)+c(\partial_i)\sigma_j\Psi+\Psi
c(\partial_i)\sigma_j]+\frac{1}{4}s+\Psi^2.\eqno(2.11)$$ By
(2.11),(2.3) and (2.4), we have
$$\omega_i=\sigma_i-\frac{1}{2}[c(\partial_i)\Psi+\Psi c(\partial_i)],\eqno(2.12)$$
$$E=-c(\partial_i)\partial^i(\Psi)-c(\partial_i)\sigma^i\Psi-\Psi
c(\partial_i)\sigma^i-\frac{1}{4}s-\Psi^2+\frac{1}{2}\partial^j[c(\partial_j)\Psi+\Psi
c(\partial_j)]$$ $$-\frac{1}{2}\Gamma^k[c(\partial_k)\Psi+\Psi
c(\partial_k)]+\frac{1}{2}\sigma^j[c(\partial_j)\Psi+\Psi
c(\partial_j)]$$ $$+\frac{1}{2}[c(\partial_j)\Psi+\Psi
c(\partial_j)]\sigma^j-\frac{g^{ij}}{4}[c(\partial_i)\Psi+\Psi
c(\partial_i)][c(\partial_j)\Psi+\Psi c(\partial_j)].\eqno(2.13)$$
For a smooth vector field $X$ on $M$, let $c(X)$ denote the Clifford action. So
$$\nabla_X=\nabla^S_X-\frac{1}{2}[c(X)\Psi+\Psi c(X)].\eqno(2.14)$$
Since $E$ is globally defined on $M$, so we can perform computations of $E$ in
normal coordinates. Taking normal coordinates about $x_0$, then
$\sigma^i(x_0)=0,~
\partial^j[c(\partial_j)](x_0)=0,~\Gamma^k(x_0)=0~g^{ij}(x_0)=\delta_i^j,$
so that
$$E(x_0)=-\frac{1}{4}s-\Psi^2+\frac{1}{2}[\partial^j(\Psi)
c(\partial_j)-c(\partial_j)\partial^j(\Psi)]-\frac{1}{4}[c(\partial_i)\Psi+\Psi
c(\partial_i)]^2(x_0)$$
$$=-\frac{1}{4}s-\Psi^2+\frac{1}{2}[e_j(\Psi)
c(e_j)-c(e_j)e_j(\Psi)]-\frac{1}{4}[c(e_i)\Psi+\Psi c(e_i)]^2(x_0)$$
$$=-\frac{1}{4}s-\Psi^2+\frac{1}{2}[\nabla^S_{e_j}(\Psi)
c(e_j)-c(e_j)\nabla^S_{e_j}(\Psi)]-\frac{1}{4}[c(e_i)\Psi+\Psi
c(e_i)]^2(x_0).\eqno(2.15)$$
We get the following Lichnerowicz formula:\\

\noindent{\bf Proposition 2.1} {\it Let $\Psi$ be a
smooth differential form on $M$ and $D_\Psi:=D+\Psi$, then}
$$D_\Psi^2=-[g^{ij}(\nabla_{\partial_i}\nabla_{\partial_j}-\nabla_{\nabla^L_{\partial_i}{\partial_j}})]
+\frac{1}{4}s+\Psi^2-\frac{1}{2}[\nabla^S_{e_j}(\Psi)
c(e_j)-c(e_j)\nabla^S_{e_j}(\Psi)]+\frac{1}{4}[c(e_i)\Psi+\Psi
c(e_i)]^2,\eqno(2.16)$${\it where $\nabla_{\partial_i}$ is defined
by (2.14) and setting $X=\partial_i$.}\\

We see two special cases of Proposition 2.1. When $\Psi=f$ where $f$ is a smooth function on $M$, we have
$$\nabla_X=\nabla_X^S-fc(X);~~E=-\frac{1}{4}s+(n-1)f^2.\eqno(2.17)$$

\noindent{\bf Corollary 2.2} {\it When $\Psi=f$, we have}
$$D_f^2=-[g^{ij}(\nabla_{\partial_i}\nabla_{\partial_j}-\nabla_{\nabla^L_{\partial_i}{\partial_j}})]
+\frac{1}{4}s+(1-n)f^2.\eqno(2.18)$$\\

 Let $\eta=a_ie^i$ be a one-form where $a_i$ is a smooth
real function and $e^i$ be a dual orthonormal frame by parallel
transport along geodesic and $X=a_ie_i$ be the dual vector field of $\eta$. When $\Psi=\sqrt{-1}c(\eta)$, by (2.14), we have
$\nabla_Y=\nabla^S_Y+\sqrt{-1}g(X,Y)$ where $Y$ is a smooth vector field on $M$. By
$e_j(c(e_i))=0$ and $de^l(x_0)=0$ (see [BGV, Lemma 4.13]), we have\\
$$E(x_0)=-\frac{1}{4}s-|X|^2+\frac{\sqrt{-1}}{2}[e_j(a^k)c(e_k)c(e_j)-c(e_j)c(e_k)e_j(a^k)]+\frac{1}{4}[c(e_i)c(X)+c(X)c(e_i)]^2$$
$$=-\frac{1}{4}s+\frac{\sqrt{-1}}{2}e_j(a^k)[c(e_k)c(e_j)-c(e_j)c(e_k)]$$
$$=-\frac{1}{4}s+\sqrt{-1}\sum_{k\neq
j}e_j(a^k)c(e_k)c(e_j)(x_0)=-\frac{1}{4}s-\sqrt{-1}c(d\eta)(x_0).\eqno(2.19)$$\\

\noindent{\bf Corollary 2.3} {\it For a one-form $\eta$ and the Clifford action $c(\eta)$, we have}
$$(D+\sqrt{-1}c(\eta))^2=-[g^{ij}(\nabla_{\partial_i}\nabla_{\partial_j}-\nabla_{\nabla^L_{\partial_i}{\partial_j}})]
+\frac{1}{4}s+\sqrt{-1}c(d\eta).\eqno(2.20)$$\\

When $\Psi$ is a two-form, we let $\Psi=2\sum_{k<l}a_{kl}e^{k}\wedge
e^{l}=\sum a_{kl}e^{k}\wedge e^{l},$ where $a_{kl}=-a_{lk},$ and
$c(\Psi)=\sum a_{kl}c(e_k)c(e_l)$. So
$$\nabla_{e_i}=e_i+\frac{1}{4}\sum_{s,t}\omega_{st}(e_i)c(e_s)c(e_t)-\sum_{k,l\neq
i}a_{kl}c(e_k)c(e_l)c(e_i),\eqno(2.21)$$ where $\omega_{st}(e_i)$
denotes the connection coefficient. By (2.15),
$$E=-\frac{1}{4}s-[a_{kl}c(e_k)c(e_l)]^2+\frac{1}{2}\left\{e_j(a_{kl})\left[c(e_k)c(e_l)c(e_j)-c(e_j)c(e_k)c(e_l)\right]\right\}$$
$$
-\frac{1}{4}\left[a_{kl}\left(c(e_i)c(e_k)c(e_l)+c(e_k)c(e_l)c(e_i)\right)\right]^2.\eqno(2.22)$$
Let $S(TM)$ be the spinor bundle on $M$ and ${\rm dim}(S(TM))=d$ and ${\rm Tr}(A)$ denote the trace of $A$ for $A\in {\rm End}(S(TM))$. Since for $k\neq l,\widetilde{k}\neq
\widetilde{l}$
$${\rm
Tr}[c(e_k)c(e_l)c(e_{\widetilde{k}})c(e_{\widetilde{l}})]=d(-\delta_k^{\widetilde{k}}\delta_l^{\widetilde{l}}
+\delta_k^{\widetilde{l}}\delta_l^{\widetilde{k}}),\eqno(2.23)$$ we
have
$${\rm
Tr}\{[a_{kl}c(e_k)c(e_l)]^2\}=-2da_{kl}^2.\eqno(2.24)$$ Since the
trace of the product of odd Clifford elements is zero, we have
$${\rm
Tr}\left[\frac{1}{2}\left\{e_j(a_{kl})\left[c(e_k)c(e_l)c(e_j)-c(e_j)c(e_k)c(e_l)\right]\right\}\right]=0.\eqno(2.25)$$
and $${\rm
Tr}\left\{\left[a_{kl}\left(c(e_i)c(e_k)c(e_l)+c(e_k)c(e_l)c(e_i)\right)\right]^2\right\}$$
$$=a_{kl}a_{\widetilde{k}\widetilde{l}}{\rm
Tr}\left[c(e_k)c(e_l)c(e_{\widetilde{k}})c(e_{\widetilde{l}})c(e_i)^2+c(e_k)c(e_l)c(e_i)^2c(e_{\widetilde{k}})c(e_{\widetilde{l}})\right.$$
$$\left.
+c(e_i)c(e_k)c(e_l)c(e_i)c(e_{\widetilde{k}})c(e_{\widetilde{l}})+c(e_k)c(e_l)c(e_i)c(e_{\widetilde{k}})c(e_{\widetilde{l}})c(e_i)\right]$$
$$=-2nda_{kl}a_{\widetilde{k}\widetilde{l}}(-\delta_k^{\widetilde{k}}\delta_l^{\widetilde{l}}+\delta_k^{\widetilde{l}}\delta_l^{\widetilde{k}})
-2\sum_{i\neq k,l}a_{kl}a_{\widetilde{k}\widetilde{l}}{\rm Tr}\left[
c(e_k)c(e_l)c(e_{\widetilde{k}})c(e_{\widetilde{l}})\right]$$
$$+2\sum_{i= k}a_{kl}a_{\widetilde{k}\widetilde{l}}{\rm Tr}\left[
c(e_k)c(e_l)c(e_{\widetilde{k}})c(e_{\widetilde{l}})\right]+
2\sum_{i= l}a_{kl}a_{\widetilde{k}\widetilde{l}}{\rm Tr}\left[
c(e_k)c(e_l)c(e_{\widetilde{k}})c(e_{\widetilde{l}})\right]$$
$$=4nda_{kl}^2+(n-2)4da_{kl}^2-4da_{kl}^2-4da_{kl}^2=8(n-2)da_{kl}^2.\eqno(2.26)$$
By (2.22) and (2.24)-(2.26), we have
$${\rm Tr}E=d\left(-\frac{1}{4}s+(6-2n)|\Psi|^2\right)\eqno(2.27)$$
and\\

\noindent{\bf Corollary 2.4} {\it let $\Psi=\sum a_{kl}e^{k}\wedge
e^{l}$ and $a_{kl}=-a_{lk}$,then ${\rm
tr}E=d\left(-\frac{1}{4}s+(6-2n)|\Psi|^2\right).$}\\

 For a general differential form $\Psi$, by (2.15) and ${\rm Tr}(AB)={\rm Tr}(BA)$, we have
$${\rm Tr}(E)={\rm Tr}\left[-\frac{1}{4}s
-\Psi^2-\frac{1}{4}[c(e_i)\Psi+\Psi c(e_i)]^2\right]$$
$$={\rm Tr}\left[-\frac{1}{4}s
-\frac{1}{2}\Psi c(e_i)\Psi
c(e_i)+(\frac{n}{2}-1)\Psi^2\right].\eqno(2.28)$$ By the
Kastler-Kalau-Walze theorem (see [Ka],[KW]), we have
$${\rm
Wres}(D_\Psi^{-n+2})=\frac{(2\pi)^{\frac{n}{2}}}{(\frac{n}{2}-2)!}\int_M{\rm
Tr}\left[\frac{1}{6}s+E\right]{\rm dvol}_M,\eqno(2.29)$$ where Wres denotes the
noncommutative residue (see [Wo]). By (2.28) and (2.29), we have\\

\noindent{\bf Theorem 2.5} {\it For even $n$-dimensional compact spin
manifolds without boundary and a general form $\Psi$, the following
equality holds:}
$${\rm
Wres}(D_\Psi^{-n+2})=\frac{(2\pi)^{\frac{n}{2}}}{(\frac{n}{2}-2)!}\int_M{\rm
Tr}\left[-\frac{1}{12}s-\frac{1}{2}c(\Psi) c(e_i)c(\Psi)
c(e_i)+(\frac{n}{2}-1)c(\Psi)^2\right]{\rm dvol}_M.\eqno(2.30)$$\\

\indent By Corollary 2.2, we have\\

\noindent {\bf Corollary 2.6} {\it For even $n$-dimensional compact spin
manifolds without boundary and a smooth function $f$ on $M$, the following equality holds:}
$${\rm
Wres}(D_f^{-n+2})=\frac{(2\pi)^{\frac{n}{2}}d}{(\frac{n}{2}-2)!}\int_M\left[-\frac{1}{12}s+(n-1)f^2\right]{\rm dvol}_M.\eqno(2.31)$$

\indent By Corollary 2.3, we have\\

\noindent {\bf Corollary 2.7} {\it For even $n$-dimensional compact spin
manifolds without boundary and a one-form $\Psi$, the following
equality holds:}
$${\rm
Wres}(D_\Psi^{-n+2})=-\frac{(2\pi)^{\frac{n}{2}}d}{12\times
(\frac{n}{2}-2)!}\int_M s{\rm dvol}_M.\eqno(2.32)$$

\indent By Corollary 2.4 and (2.29), we have\\

\noindent {\bf Corollary 2.8} {\it For even $n$-dimensional compact spin
manifolds without boundary and a two-form $\Psi$, the following
equality holds:}
$${\rm
Wres}(D_\Psi^{-n+2})=\frac{(2\pi)^{\frac{n}{2}}d}{(\frac{n}{2}-2)!}\int_M{\rm
Tr}\left[-\frac{1}{12}s+(6-2n)|\Psi|^2\right]{\rm
dvol}_M.\eqno(2.33)$$\\

\noindent{\bf 2.2 A Kastler-Kalau-Walze type theorem for
perturbations of Dirac operators on manifolds with boundary}\\

\indent  We now let $M$ be a compact $4$-dimensional spin manifold with boundary $\partial M$ and $U\subset M$ be the collar neighborhood of $\partial M$ which is diffeomorphic to
$\partial M\times [0,1).$
and we shall compute the noncommutative residue for manifolds with boundary of $(\pi^+D_\Psi^{-1})^2$. That is, we shall compute
 $\widetilde{{\rm Wres}}[(\pi^+D_\Psi^{-1})^2]$ (for the related definitions, see
 [Wa1]) and we take the metric as in [Wa1]. Let $(x',x_n)\in U$ where $x'\in \partial M$ and $x_n$ denotes the normal direction coordinate.
By (2.2.4) in [Wa1], we have
$$\widetilde{{\rm
Wres}}[(\pi^+D_\Psi^{-1})^2]=\int_M\int_{|\xi|=1}{\rm
trace}_{S(TM)}[\sigma_{-4}(D_\Psi^{-2})]\sigma(\xi)dx+\int_{\partial
M}\Phi,\eqno(2.34)$$ \indent where
$$\Phi=\int_{|\xi'|=1}\int^{+\infty}_{-\infty}\sum^{\infty}_{j, k=0}
\sum\frac{(-i)^{|\alpha|+j+k+1}}{\alpha!(j+k+1)!}
{\rm trace}_{S(TM)}
[\partial^j_{x_n}\partial^\alpha_{\xi'}\partial^k_{\xi_n}
\sigma^+_{r}(D_\Psi^{-1})(x',0,\xi',\xi_n)$$
$$\times
\partial^\alpha_{x'}\partial^{j+1}_{\xi_n}\partial^k_{x_n}\sigma_{l}
(D_\Psi^{-1})(x',0,\xi',\xi_n)]d\xi_n\sigma(\xi')dx',\eqno(2.35)$$
\noindent where the sum is taken over $
r-k-|\alpha|+l-j-1=-4,~~r,l\leq-1$ and
$\sigma^+_{r}(D_\Psi^{-1})=\pi^+_{\xi_n}\sigma_{r}(D_\Psi^{-1})$
(for the definition of  $\pi^+$, see [Wa1]). By Theorem 2.5, we have
$$\int_M\int_{|\xi|=1}{\rm tr}[\sigma_{-4}(D_\Psi^{-2})]\sigma(\xi)dx={4\pi^2}\int_M{\rm
Tr}\left[-\frac{1}{12}s-\frac{1}{2}\Psi c(e_i)\Psi
c(e_i)+(\frac{n}{2}-1)\Psi^2\right]{\rm dvol}_M.\eqno(2.36)$$ So we
only need to compute $\int_{\partial M}\Phi$. In analogy with Lemma 2.1 of
[Wa1], we can prove the following useful result.\\

\noindent{\bf Lemma 2.9} {\it The symbolic calculus of pseudodifferential operators yields}
$$
q_{-1}(D_\Psi^{-1})=\frac{\sqrt{-1}c(\xi)}{|\xi|^2},~~q_{-2}(D_\Psi^{-1})=q_{-2}(D^{-1})+\frac{c(\xi)\Psi
c(\xi)}{|\xi|^4}.\eqno(2.37)$$\\

Similar to the computations in section 2.2.2 in [Wa1], we can
split $\Phi$ into the sum of five terms. Since
$q_{-1}(D_\Psi^{-1})=q_{-1}(D^{-1})$, then terms (a)(I), (II),(III)
in our case are the same as the terms (a)(I), (II),(III) in [Wa1], so
$${\rm term~(a)~I)}=-\int_{|\xi'|=1}\int^{+\infty}_{-\infty}\sum_{|\alpha|=1}
{\rm trace} [\partial^\alpha_{\xi'}\pi^+_{\xi_n}q_{-1}\times
\partial^\alpha_{x'}\partial_{\xi_n}q_{-1}](x_0)d\xi_n\sigma(\xi')dx'=0.\eqno(2.38)$$
$${\rm term~
(a)~II)}=-\frac{1}{2}\int_{|\xi'|=1}\int^{+\infty}_{-\infty} {\rm
trace} [\partial_{x_n}\pi^+_{\xi_n}q_{-1}\times
\partial_{\xi_n}^2q_{-1}](x_0)d\xi_n\sigma(\xi')dx'\eqno(2.39)$$
$~~~~~~~~~~~~~~~~~~~~~=-\frac{3}{8}\pi h'(0)\Omega_3dx'.$

$${\rm term~ (a)~III)}=-\frac{1}{2}\int_{|\xi'|=1}\int^{+\infty}_{-\infty}
{\rm trace} [\partial_{\xi_n}\pi^+_{\xi_n}q_{-1}\times
\partial_{\xi_n}\partial_{x_n}q_{-1}](x_0)d\xi_n\sigma(\xi')dx'\eqno(2.40)$$
$~~~~~~~~~~~~~~~~~~~~=\frac{3}{8}\pi h'(0)\Omega_3dx'.$\\
Then we only need to compute the term (b) and the term (c). By Lemma
2.9,\\
$${\rm term~ (b)}:=-i\int_{|\xi'|=1}\int^{+\infty}_{-\infty}
{\rm trace} [\pi^+_{\xi_n}q_{-2}(D^{-1})\times
\partial_{\xi_n}q_{-1}(D^{-1})](x_0)d\xi_n\sigma(\xi')dx'$$
$$
-i\int_{|\xi'|=1}\int^{+\infty}_{-\infty} {\rm trace}
[\pi^+_{\xi_n}\left(\frac{c(\xi)\Psi c(\xi)}{|\xi|^4}\right)\times
\partial_{\xi_n}q_{-1}(D^{-1})](x_0)d\xi_n\sigma(\xi')dx'.
\eqno(2.41)$$ By the term (b) in [Wa1], we have
$$-i\int_{|\xi'|=1}\int^{+\infty}_{-\infty}
{\rm trace} [\pi^+_{\xi_n}q_{-2}(D^{-1})\times
\partial_{\xi_n}q_{-1}(D^{-1})](x_0)d\xi_n\sigma(\xi')dx'=\frac{9}{8}\pi
h'(0)\Omega_3dx',\eqno(2.42)$$ where $\Omega_3$ is the canonical
volume of $3$-dimensional unit sphere. Moreover
$$
\pi^+_{\xi_n}\left(\frac{c(\xi)\Psi
c(\xi)}{|\xi|^4}\right)(x_0)|_{|\xi'|=1}
=\pi^+_{\xi_n}\left[\frac{[c(\xi')+\xi_nc(dx_n)]\Psi[c(\xi')+\xi_nc(dx_n)]}{(1+\xi_n^2)^2}\right]$$
$$=\frac{1}{2\pi i}\int_{\Gamma^+}\frac{\frac{c(\xi')\Psi
c(\xi')+c(dx_n)\Psi c(\xi')\eta_n+ c(\xi')\Psi
c(dx_n)\eta_n+c(dx_n)\Psi c(dx_n)\eta_n^2
}{(\eta_n+i)^2(\xi_n-\eta_n)}} {(\eta_n-i)^2}d\eta_n$$
$$=\left[\frac{c(\xi')\Psi c(\xi')+c(dx_n)\Psi c(\xi')\eta_n+
c(\xi')\Psi c(dx_n)\eta_n+c(dx_n)\Psi c(dx_n)\eta_n^2
}{(\eta_n+i)^2(\xi_n-\eta_n)} \right]^{(1)}|_{\eta_n=i}$$
$$
=-\frac{i\xi_n+2}{4(\xi_n-i)^2}c(\xi')\Psi
c(\xi')-\frac{i}{4(\xi_n-i)^2}[c(dx_n)\Psi c(\xi')$$ $$+c(\xi')\Psi
c(dx_n)]-\frac{i\xi_n}{4(\xi_n-i)^2}c(dx_n)\Psi
c(dx_n).\eqno(2.43)$$

$$\partial_{\xi_n}q_{-1}|_{|\xi'|=1}=\sqrt{-1}\left[\frac{1-\xi_n^2}{(1+\xi_n^2)^2}c(dx_n)-\frac{2\xi_n}{(1+\xi_n^2)^2}c(\xi')\right].\eqno(2.44)$$

\noindent By (2.43) and (2.44) and
$$c(\xi')^2|_{|\xi'|=1}=-1,~~c(dx_n)^2=-1,
~c(\xi')c(dx_n)=-c(dx_n)c(\xi'), {\rm Tr}(AB)={\rm
Tr}(BA),\eqno(2.45)$$ we get
$$ {\rm trace}
\left[\pi^+_{\xi_n}\left(\frac{c(\xi)\Psi c(\xi)}{|\xi|^4}\right)\times
\partial_{\xi_n}q_{-1}(D^{-1})\right](x_0)|_{|\xi'|=1}$$
$$=
\frac{\sqrt{-1}}{2(1+\xi_n^2)^2}{\rm
Tr}[c(dx_n)\Psi]+\frac{1}{2(1+\xi_n^2)^2}{\rm
Tr}[c(\xi')\Psi].\eqno(2.46)$$
Considering for $i<n$, $\int_{|\xi'|=1}\xi_i\sigma(\xi')=0$, then\\
$$
-i\int_{|\xi'|=1}\int^{+\infty}_{-\infty} {\rm trace}
[\pi^+_{\xi_n}\left(\frac{c(\xi)\Psi c(\xi)}{|\xi|^4}\right)\times
\partial_{\xi_n}q_{-1}(D^{-1})](x_0)d\xi_n\sigma(\xi')dx'$$
$$=\frac{\pi}{4}\Omega_3{\rm Tr}[c(dx_n)\Psi]dx',
\eqno(2.47)$$ and
$${\rm term~ (b)}=\frac{9}{8}\pi
h'(0)\Omega_3dx'+\frac{\pi}{4}\Omega_3{\rm Tr}[c(dx_n)\Psi]dx'.\eqno(2.48)$$\\
Similarly, we have
$${\rm term~ (c)}=-\frac{9}{8}\pi
h'(0)\Omega_3dx'-\frac{\pi}{4}\Omega_3{\rm Tr}[c(dx_n)\Psi]dx'.\eqno(2.49)$$\\
Then the sum of terms (b) and (c) is zero and $\Phi$ is zero. Then we get\\

\noindent {\bf Theorem 2.10}~~{\it  Let $M$ be a $4$-dimensional
compact spin manifold with boundary $\partial M$ and the metric
$g^M$ (see (1.3) in [Wa1]). Let $\Psi$ be a general differential form on $M$. Then}
$$ \widetilde{{\rm
Wres}}[(\pi^+D_\Psi^{-1})^2]={4\pi^2}\int_M{\rm
Tr}\left[-\frac{1}{12}s-\frac{1}{2}c(\Psi) c(e_i)c(\Psi)
c(e_i)+c(\Psi)^2\right]{\rm dvol}_M.\eqno(2.50)$$\\

\indent In [Wa2], we proved a Kastler-Kalau-Walze theorem
associated to Dirac operators for $6$-dimensional spin manifolds
with boundary.
In fact, our computations hold for general Laplacians. This implies\\

\noindent {\bf Proposition 2.11 ([Wa2])}~~{\it Let $M$ be a
$6$-dimensional compact Riemannian manifold with boundary
$\partial M$ and the metric $g^M$ (see (1.3) in [Wa1]). Let $\Delta$ be a general
Laplacian acting on sections of the vector bundle $V$. Then}
$$ \widetilde{{\rm
Wres}}[(\pi^+\Delta^{-1})^2]=8\pi^3\int_M{\rm Tr}[\frac{s}{6}+E]{\rm
dvol}_M.\eqno(2.51)$$

Since $D_\Psi^2$ is a general Laplacian, then we get \\

\noindent {\bf Corollary 2.12}~~{\it Let $M$ be a $6$-dimensional
compact spin manifold with boundary $\partial M$ and the metric
$g^M$. Let $\Psi$ be a general differential form on $M$. Then}
$$ \widetilde{{\rm
Wres}}[(\pi^+D_\Psi^{-2})^2]=8\pi^3\int_M{\rm
Tr}\left[-\frac{1}{12}s-\frac{1}{2}\Psi c(e_i)\Psi
c(e_i)+2\Psi^2\right]{\rm dvol}_M.\eqno(2.52)$$\\

In the above two cases, the boundary terms vanish. In the following, we will give a boundary term nonvanishing case and compute
${\rm
Wres}((D_\Psi D)^{-1}).$ We have\\
$$D_\Psi D=-g^{ij}\partial_i\partial_j+\left(-2\sigma^j+\Gamma^j+\Psi
c(\partial^j)\right)\partial_j$$
$$+g^{ij}[-\partial_i(\sigma_j)-\sigma_i\sigma_j+\Gamma_{ij}^k\sigma_k
+\Psi c(\partial_i)\sigma_j]+\frac{1}{4}s,\eqno(2.53)$$ and
$$\omega_i=\sigma_i-\frac{1}{2}\Psi c(\partial_i),\eqno(2.54)$$
$$E=-\Psi
c(\partial_i)\sigma^i-\frac{1}{4}s+\frac{1}{2}\partial^j[\Psi
c(\partial_j)]$$ $$-\frac{1}{2}\Gamma^k\Psi
c(\partial_k)+g^{ij}\left[\frac{1}{2}\sigma_i\Psi
c(\partial_j)+\frac{1}{2}\Psi c(\partial_i)\sigma_j -\frac{1}{4}\Psi
c(\partial_i)\Psi c(\partial_j)\right].\eqno(2.55)$$
Similar to the proof of (2.15), we have
$$E=-\frac{1}{4}s+\frac{1}{2}\nabla^S_{e_i}(\Psi)c(e_i)-\frac{1}{4}\Psi c(e_i)\Psi c(e_i).\eqno(2.56)$$
So
$$D_\Psi D=-[g^{ij}(\nabla_{\partial_i}\nabla_{\partial_j}-\nabla_{\nabla^L_{\partial_i}{\partial_j}})]
+\frac{1}{4}s-\frac{1}{2}\nabla^S_{e_i}(\Psi)c(e_i)+\frac{1}{4}\Psi
c(e_i)\Psi c(e_i).\eqno(2.57)$$
Then we get\\

\noindent {\bf Proposition 2.13}~~{\it  Let $M$ be a $4$-dimensional
compact spin manifold without boundary. Then}
$$ {\rm
Wres}[(D_\Psi D)^{-1}]={4\pi^2}\int_M{\rm
Tr}\left[-\frac{1}{12}s+\frac{1}{2}\nabla^S_{e_i}(\Psi)c(e_i)-\frac{1}{4}\Psi
c(e_i)\Psi c(e_i)
\right]{\rm dvol}_M.\eqno(2.58)$$\\

\indent When $\Psi$ is a one-form, we can get the following
 corollary:\\

\noindent {\bf Corollary 2.14}~~{\it  Let $M$ be a $4$-dimensional
compact spin manifold without boundary and  let $\Psi$ be a one-form on $M$. Then}
$$ {\rm
Wres}[(D_\Psi D)^{-1}]={16\pi^2}\int_M\left[-\frac{1}{12}s+\frac{1}{2}\delta(\Psi)-2|\Psi|^2\right]{\rm dvol}_M.\eqno(2.59)$$\\

Now we compute $\widetilde{{\rm
Wres}}[\pi^+D_\Psi^{-1}\pi^+D^{-1}]$. We have that terms (a) and (b) are
the same as in Theorem 2.10, and since ${\rm term~ (c)}=-\frac{9}{8}\pi
h'(0)\Omega_3dx'$, we get
$$\int_{\partial M}\Phi=\frac{\pi}{4}\Omega_3\int_{\partial M}{\rm Tr}[c(dx_n)\Psi]{\rm dvol}_{\partial_M},\eqno(2.60)$$
and\\

\noindent {\bf Proposition 2.15}~~{\it  Let $M$ be a $4$-dimensional
compact spin manifold with boundary. Then}
$$ \widetilde{{\rm
Wres}}[\pi^+D_\Psi^{-1}\pi^+D^{-1}]={4\pi^2}\int_M{\rm
Tr}\left[-\frac{1}{12}s+\frac{1}{2}e_i(\Psi)c(e_i)-\frac{1}{4}\Psi c(e_i)\Psi c(e_i)
\right]{\rm dvol}_M$$
$$+\frac{\pi}{4}\Omega_3\int_{\partial M}{\rm Tr}[c(dx_n)\Psi]{\rm dvol}_{\partial_M}.\eqno(2.61)$$\\

\noindent {\bf Remark.} When $\Psi$ is not a one-form,
then the boundary term vanishes. When $\Psi=Kdx_n$ near the boundary
where $K$ is the extrinsic curvature, then the boundary term is proportional
 to the gravitational action on the boundary. In fact, the reason of the boundary term
 being not zero is that $\pi^+D_\Psi$ and $\pi^+D$ are not symmetric.\\

\section{A Kastler-Kalau-Walze type theorem for conformal perturbations of Dirac operators }

   \quad In [CM], Connes and Moscovici defined a twisted spectral triple and considered the conformal
    Dirac operator $e^{h}De^{h}$ where $h$ is a smooth function on a manifold $M$ without boundary. We want to compute ${\rm
    Wres}[(e^{h}De^{h})^{-2}].$ We know that
       $${\rm
    Wres}[(e^{h}De^{h})^{-2}]={\rm
    Wres}[e^{-h}D^{-1}e^{-2h}D^{-1}e^{-h}]={\rm
    Wres}[e^{-2h}D^{-1}e^{-2h}D^{-1}].\eqno(3.1)$$
In the following, we will compute the more general case, i.e.
   ${\rm Wres}[fD^{-1}gD^{-1}]$ for nonzero smooth functions $f,g$
   and prove a Kastler-Kalau-Walze type theorem for conformal Dirac operators. When $f=g=e^{-2h}$, we get the expression of ${\rm
    Wres}[(e^{h}De^{h})^{-2}].$ We have

   \begin{eqnarray*}
   {\rm Wres}[fD^{-1}gD^{-1}]&=&{\rm Wres}[(f^{-1}Dg^{-1}D)^{-1}]\\
   &=&{\rm Wres}\left\{(f^{-1}g^{-1}D^2+f^{-1}[D,g^{-1}]D)^{-1}\right\}\\
&=&\int_Mfg{\rm wres}[(D^2-g^{-1}c(dg)D)^{-1}],~~~~~~~~~~~~~~~~~~~~~~~~~~~~~~~(3.2)
\end{eqnarray*}
where ${\rm wres}$ denotes the residue density and we note that the Kastler-Kalau-Walze theorem holds at the residue density level. Some computations show that \\
\begin{eqnarray*}
D^2-g^{-1}c(dg)D&=&-g^{ij}\partial_i\partial_j+[-2\sigma^j+\Gamma^j-g^{-1}c(dg)c(\partial^j)]\partial_j\\
&&+[-\partial^j\sigma_j-\sigma^j\sigma_j+
\Gamma^k\sigma_k+\frac{1}{4}s-g^{-1}c(dg)c(\partial^j)\sigma_j].~~~~~~~~~~(3.3)\\
\omega_i&=&\sigma_i+\frac{1}{2}g^{-1}c(dg)c(\partial_i),~~~~~~~~~~~~~~~~~~~~~~~~~~~~~~~~~~~~~~~~~~~~~~~(3.4)\\
E&=&-\frac{s}{4}+g^{-1}c(dg)c(\partial^j)\sigma_j-\partial^j[\frac{1}{2}g^{-1}c(dg)c(\partial_j)]\\
&&-\frac{1}{2}\sigma^jg^{-1}c(dg)c(\partial_j)
-\frac{1}{2}g^{-1}c(dg)c(\partial_i)\sigma^i\\
&&-\frac{1}{4}g^{ij}g^{-1}c(dg)c(\partial_i)g^{-1}c(dg)c(\partial_j)+\frac{1}{2}g^{-1}c(dg)c(\partial_k)
\Gamma^k.~~~(3.5)
\end{eqnarray*}\\

 Since $E$ is globally defined, we can compute it
in the normal coordinates. Then we have
$${\rm Tr}(E)(x_0)={\rm Tr}\left[-\frac{s}{4}-\frac{1}{2}\partial_j(g^{-1}c(dg))c(\partial_j)
-\frac{1}{4}g^{-1}c(dg)c(\partial_i)g^{-1}c(dg)c(\partial_i)\right](x_0),\eqno(3.6)$$
and
\begin{eqnarray*}
&&
{\rm Tr}\left[g^{-1}c(dg)c(\partial_i)g^{-1}c(dg)c(\partial_i)\right](x_0)\\
&=&g^{-2}{\rm Tr}\left[\sum_{i,k,l}\frac{\partial g}{\partial x_k}\frac{\partial g}{\partial x_l}
c(\partial_k)c(\partial_i)c(\partial_l)c(\partial_i)\right]\\
&=&g^{-2}{\rm Tr}\left[\sum_{i,k}\left(\frac{\partial g}{\partial x_k}\right)^2
c(\partial_k)c(\partial_i)c(\partial_k)c(\partial_i)\right]\\
&=&g^{-2}{\rm Tr}\left[\sum_{i\neq k}\left(\frac{\partial g}{\partial x_k}\right)^2
c(\partial_k)c(\partial_i)c(\partial_k)c(\partial_i)+\sum_k\left(\frac{\partial g}{\partial x_k}\right)^2c(\partial_k)^4\right]\\
&=&-2g^{-2}\sum_k\left(\frac{\partial g}{\partial x_k}\right)^2{\rm Tr}[{\rm
Id}]. ~~~~~~~~~~~~~~~~~~~~~~~~~~~~~~~~~~~~~~~~~~~~~~~~~~~~~~~~~~~~~~(3.7)
\end{eqnarray*}
 Similarly,
$${\rm Tr}\left[\partial_j(g^{-1}c(dg))c(\partial_j)\right]=\sum_j\left[\frac{1}{g^2}\left(\frac{\partial g}{\partial x_j}\right)^2
-g^{-1}\frac{\partial^2 g}{\partial x_j^2}\right].\eqno(3.8)$$ So
$${\rm Tr}[\frac{s}{6}+E]=-\frac{s}{3}+2g^{-1}\sum_j\frac{\partial^2 g}{\partial x_j^2}=-\frac{s}{3}-2g^{-1}\triangle(g).\eqno(3.9)$$
  By $$\int_Mf\triangle(g){\rm dvol}_M=\int_M\left<df,dg\right>{\rm dvol}_M,\eqno(3.10)$$
 we get\\

\noindent {\bf Theorem 3.1}~~{\it  Let $M$ be a $4$-dimensional
compact spin manifold without boundary, then}
$${\rm Wres}[fD^{-1}gD^{-1}]=-{4\pi^2}\int_M[\frac{fgs}{3}+2\left<df,dg\right>]
{\rm dvol}_M.\eqno(3.11)$$\\

\noindent{\bf Remark.} In Theorem 3.1, when $f=g=e^{-2h}$, we get a
Kastler-Kalau-Walze theorem for conformal Dirac operators. In fact, Theorem 3.1 holds
true for any choice of the smooth functions $f$ and $g$, since we can prove
(3.11) by means of the symbolic calculus of pseudodifferential operators
without using (3.2), and it is not essential that $f$ and $g$ are nowhere
vanishing.\\

 Now we consider
manifolds with boundary and we will compute\\ $ \widetilde{{\rm
Wres}}[\pi^+(fD^{-1})\pi^+(gD^{-1})].$ As in [Wa1], we have five
terms.
$${\rm term~ (a)~ I}=-fg\int_{|\xi'|=1}\int^{+\infty}_{-\infty}\sum_{|\alpha|=1}
{\rm trace} [\partial^\alpha_{\xi'}\pi^+_{\xi_n}q_{-1}\times
\partial^\alpha_{x'}\partial_{\xi_n}q_{-1}](x_0)d\xi_n\sigma(\xi')dx'$$
$$-f\sum_{j<n}\partial_j(g)\int_{|\xi'|=1}\int^{+\infty}_{-\infty}\sum_{|\alpha|=1}
{\rm trace} [\partial_{\xi_j}\pi^+_{\xi_n}q_{-1}\times
\partial_{\xi_n}q_{-1}](x_0)d\xi_n\sigma(\xi')dx'=0,\eqno(3.12)$$
$${\rm term~ (a)~ II}=-\frac{1}{2}fg\int_{|\xi'|=1}\int^{+\infty}_{-\infty}
{\rm trace} [\partial_{x_n}\pi^+_{\xi_n}q_{-1}\times
\partial_{\xi_n}^2q_{-1}](x_0)d\xi_n\sigma(\xi')dx'$$
$$-\frac{1}{2}g\partial_{x_n}f\int_{|\xi'|=1}\int^{+\infty}_{-\infty}
{\rm trace} [\pi^+_{\xi_n}q_{-1}\times
\partial_{\xi_n}^2q_{-1}](x_0)d\xi_n\sigma(\xi')dx'$$
$$=-\frac{3}{8}\pi h'(0)\Omega_3fgdx'-\frac{\pi
i}{2}\Omega_3g\partial_{x_n}(f)dx'.\eqno(3.13)$$
$${\rm term ~(a)~ III)}=\frac{3}{8}\pi h'(0)\Omega_3fgdx'+\frac{\pi i}{2}\Omega_3f\partial_{x_n}(g)dx'.\eqno(3.14)$$
As in [Wa1], we have
$${\rm term~ (b)}=-i\int_{|\xi'|=1}\int^{+\infty}_{-\infty}
{\rm trace} [\pi^+_{\xi_n}q_{-2}\times
\partial_{\xi_n}q_{-1}](x_0)d\xi_n\sigma(\xi')dx'
=\frac{9}{8}fg\pi h'(0)\Omega_3dx'.\eqno(3.15)$$
and
$${\rm term~ (c)}=-i\int_{|\xi'|=1}\int^{+\infty}_{-\infty}
{\rm trace} [\pi^+_{\xi_n}q_{-1}\times
\partial_{\xi_n}q_{-2}](x_0)d\xi_n\sigma(\xi')dx'
=-\frac{9}{8}fg\pi h'(0)\Omega_3dx'.\eqno(3.16)$$
So the sum of terms (b) and (c) is zero. Then we obtain
$$\int_{\partial M}\Phi=\frac{\pi i\Omega_3}{2}\int_{\partial M}[f\partial_{x_n}(g)-g\partial_{x_n}(f)]_{x_n=0}{\rm dvol}_{\partial M}.\eqno(3.17)$$\\
By the definition of the noncommutative residue for manifolds with boundary, we have that the interior term of $\widetilde{{\rm
Wres}}[\pi^+(fD^{-1})\pi^+(gD^{-1})] $ equals $ {\rm Wres}[fD^{-1}gD^{-1}]$. Then by Theorem 3.1 and (3.17), we get\\

\noindent {\bf Theorem 3.2}~~{\it  Let $M$ be a $4$-dimensional
compact spin manifold with boundary. Then}
$$\widetilde{{\rm
Wres}}[\pi^+(fD^{-1})\pi^+(gD^{-1})]=-{4\pi^2}\int_M[\frac{fgs}{3}+2\left<df,dg\right>]
{\rm dvol}_M$$\\
$$+\frac{\pi i\Omega_3}{2}\int_{\partial M}[f\partial_{x_n}(g)-g\partial_{x_n}(f)]|_{x_n=0}{\rm dvol}_{\partial M}.\eqno(3.18)$$\\

\indent When $f=1,~g=x_nK$ near the boundary, we have that the boundary term is proportional to the gravitational action on the boundary.\\

\section{The spectral action for perturbations of Dirac operators}

  \quad In [IL], Iochum and Levy computed heat kernel coefficients for Dirac operators with one form perturbations
     and proved that there are no tadpoles for compact spin manifolds without boundary. In [SZ], they investigated the spectral action for scalar perturbations of Dirac operators. In [HPS], Hanisch, Pf\"{a}ffle and
    Stephan derived a formula for the gravitational part of the
    spectral action for Dirac operators on $4$-dimensional spin
    manifolds with totally anti-symmetric torsion. In fact Dirac operators with totally anti-symmetric torsion are three form perturbations of Dirac operators. In this section, we will give some details on the spectral action for Dirac operators with scalar perturbations. We also compute the spectral action for Dirac operators with two-form perturbations on $4$-dimensional compact spin
      manifolds without boundary.\\
  \indent For the perturbed self-adjoint Dirac operator $D_{\Psi}$, we will calculate
the bosonic part of the spectral action. It is defined to be the
number of eigenvalues of $D_{\Psi}$ in the interval
$[-\wedge,\wedge]$ with $\wedge\in {\bf R}^+$. As in [CC1], it is
expressed as
$$ I={\rm Tr}{F}\left(\frac{D^2_\Psi}{\wedge^2}\right).\eqno(4.1)$$
\noindent Here Tr denotes the operator trace in the $L^2$ completion
of $\Gamma (M,S(TM))$, and ${F}:{\bf R}^+\rightarrow {\bf R}^+$ is a
cut-off function with support in the interval $[0,1]$ which is
constant near the origin. Let ${\rm dim}~M=n$. By Lemma 1.7.4 in [Gi], we
have the heat trace asymptotics for $t\rightarrow 0$,
$${\rm Tr}(e^{-tD_\Psi^2})\sim \sum_{m\geq
0}t^{m-\frac{n}{2}}a_{2m}(D_\Psi^2).\eqno(4.2)$$ One uses the
Seeley-deWitt coefficients $a_{2m}(D_\Psi^2)$ and $t=\wedge^{-2}$ to
obtain an asymptotics for the spectral action when ${\rm dim}~ M=4$
[CC1]
$$I={\rm tr}{F}\left(\frac{D^2_\Psi}{\wedge^2}\right)\sim
\wedge^4F_4a_0(D^2_\Psi)+\wedge^2F_2a_2(D^2_\Psi)+\wedge^0F_0a_4(D^2_\Psi)~~{\rm
as} ~~\wedge\rightarrow \infty \eqno(4.3)$$ \noindent with the first
three moments of the cut-off function which are given by
$F_4=\int_0^{\infty}s{F}(s)ds,$\\
$F_2=\int_0^{\infty}{F}(s)ds$ and $F_0={F}(0)$.
 Let
$$\Omega_{ij}={\nabla}_{e_i}{\nabla}_{e_j}-{\nabla}_{e_j}{\nabla}_{e_i}
-{\nabla}_{[e_i,e_j]}.\eqno(4.4)$$ We use [Gi, Thm 4.1.6]
to obtain the first three coefficients of the heat trace
asymptotics:
$$a_0(D_\Psi)=(4\pi)^{-\frac{n}{2}}\int_M{\rm Tr}(\rm
Id)dvol,\eqno(4.5)$$
$$a_2(D_\Psi)=(4\pi)^{-\frac{n}{2}}\int_M{\rm
Tr}[\frac{s}{6}+E]dvol,\eqno(4.6)$$
$$a_4(D_\Psi)=\frac{(4\pi)^{-\frac{n}{2}}}{360}\int_M{\rm
Tr}[-12R_{ijij,kk}+5R_{ijij}R_{klkl}$$
$$-2R_{ijik}R_{ljlk}+2R_{ijkl}R_{ijkl}-60R_{ijij}E+180E^2+60E_{,kk}+30\Omega_{ij}
\Omega_{ij}]dvol.\eqno(4.7)$$\\
 \indent When $\Psi=f$, by (2.17) and (4.6),
 $$a_2(D_f)=(2\pi)^{-\frac{n}{2}}[-\frac{s}{12}+(n-1)f^2].\eqno(4.8)$$
 $$5s^2+60sE+180E^2=\frac{5}{4}s^2-30(n-1)sf^2+180(n-1)^2f^4.\eqno(4.9)$$
${\rm Tr}[\Omega_{ij}\Omega_{ij}]$ is globally defined, thus we only
compute it in normal coordinates about $x_0$ and the local
orthonormal frame $e_i$ obtained by parallel tansport along
geodesics from $x_0$. Then
$$\omega_{st}(x_0)=0,~\partial_i(c(e_j))=0,~[e_i,e_j](x_0)=0.\eqno(4.10)$$
We know that the curvature of the canonical spin connection is
$$R^S(e_i,e_j)=-\frac{1}{4}\sum_{s,t=1}^nR^M_{ijst}c(e_s)c(e_t).$$ Then we have
$$\Omega(e_i,e_j)(x_0)=[e_i+\frac{1}{4}\sum_{s,t}\omega_{st}(e_i)c(e_s)c(e_t)-fc(e_i)][e_j+\frac{1}{4}\sum_{s,t}\omega_{st}(e_j)c(e_s)c(e_t)-fc(e_j)]$$
$$-[e_j+\frac{1}{4}\sum_{s,t}\omega_{st}(e_j)c(e_s)c(e_t)-fc(e_i)][e_i+\frac{1}{4}\sum_{s,t}\omega_{st}(e_i)c(e_s)c(e_t)-fc(e_i)]$$
$$=-\frac{1}{4}\sum_{s,t=1}^nR^M_{ijst}c(e_s)c(e_t)-e_i(f)c(e_j)+e_j(f)c(e_i)+2f^2c(e_i)c(e_j),~
{\rm for}~ i\neq j.\eqno(4.11)$$ So
$${\rm Tr}[\Omega_{ij}\Omega_{ij}](x_0)=\sum_{i\neq j}{\rm
Tr}\left\{\frac{1}{16}\sum_{s,t,s_1,t_1=1}^nR^M_{ijst}R^M_{ijs_1t_1}c(e_s)c(e_t)c(e_{s_1})c(e_{t_1})+e_i(f)^2c(e_j)^2\right.$$
$$+e_j(f)^2c(e_i)^2+4f^4c(e_i)c(e_j)c(e_i)c(e_j)$$
$$
-\frac{f^2}{2}\sum_{s,t=1}^nR^M_{ijst}[c(e_s)c(e_t)c(e_i)c(e_j)
+c(e_i)c(e_j)c(e_s)c(e_t)]$$
$$\left.-e_i(f)e_j(f)[c(e_j)c(e_i)+c(e_i)c(e_j)]\right\}.\eqno(4.12)$$ By
(2.23), we obtain
$$\sum_{i,j,s,t,s_1,t_1=1}^n{\rm
Tr}[\frac{1}{16}R^M_{ijst}R^M_{ijs_1t_1}c(e_s)c(e_t)c(e_{s_1})c(e_{t_1})]=-\frac{d}{8}\sum_{i,j,s,t=1}^n(R^M_{ijst})^2,\eqno(4.13)$$
$$\sum_{i\neq j}{\rm
Tr}[e_i(f)^2c(e_j)^2+e_j(f)^2c(e_i)^2]=2d(1-n)\sum_ie_i(f)^2=2d(1-n)|df|^2,\eqno(4.14)$$
$$\sum_{i\neq j}{\rm
Tr}[4f^4c(e_i)c(e_j)c(e_i)c(e_j)]=-4dn(n-1)f^4,\eqno(4.15)$$
$$\sum_{i\neq j}{\rm
Tr}[-e_i(f)e_j(f)(c(e_j)c(e_i)+c(e_i)c(e_j))]=0,\eqno(4.16)$$
$$\sum_{i\neq j}{\rm
Tr}\left\{-\frac{f^2}{2}R^M_{ijst}[c(e_s)c(e_t)c(e_i)c(e_j)
+c(e_i)c(e_j)c(e_s)c(e_t)]\right\}=-2f^2ds.\eqno(4.17)$$ By
(4.12)-(4.17), we obtain
$${\rm
Tr}[\Omega_{ij}\Omega_{ij}]=-\frac{d}{8}(R^M_{ijst})^2+2d(1-n)|df|^2-2f^2ds-4dn(n-1)f^4.\eqno(4.18)$$
By (4.7) (4.9) and (4.18), we get\\

\noindent {\bf Proposition 4.1 ([SZ])}{\it The following equality holds}
$$a_4(D_f)=\frac{d}{360\times (4\pi)^{\frac{n}{2}}}\left[3\Delta
s+\frac{5}{4}s^2-30(n+1)sf^2+60(n-1)(n-3)f^4\right.$$
$$\left.-2R_{ijik}R_{ljlk}-\frac{7}{4}R_{ijst}^2+60(1-n)|df|^2-60(n-1)\Delta(f^2)\right].\eqno(4.19)$$\\

\indent In the following, we assume that ${\rm dim}~M=4$ and $d=4$.
We let $\Psi$ be a two-form, namely $\Psi=\sum_{k,l}a_{kl}e^k\wedge e^l$
where $a_{kl}=-a_{lk}$. We may consider $\sqrt{-1}\Psi$ for
selfadjoint perturbed Dirac operators. By Corollary 2.4, we obtain
$$a_2(D_\Psi)=d(4\pi)^{-\frac{n}{2}}[-\frac{s}{12}+(6-2n)|\Psi|^2]=\frac{-1}{4\pi^2}[\frac{s}{12}+2|\Psi|^2].\eqno(4.20)$$
Firstly, we compute ${\rm Tr}(E^2)$. By (2.22) and (4.10),
$$\sum_{j,k,l} e_j(a_{kl})\left[c(e_k)c(e_l)c(e_j)-c(e_j)c(e_k)c(e_l)\right]=4e_k(a_{kl})c(e_l),\eqno(4.21)$$
\begin{eqnarray*}
&&a_{kl}\left[c(e_i)c(e_k)c(e_l)+c(e_k)c(e_l)c(e_i)\right]a_{k_1l_1}\left[c(e_i)c(e_{k_1}
)c(e_{l_1})+c(e_{k_1})c(e_{l_1})c(e_i)\right]\\
&=&\sum_{i\neq k,l}2a_{kl}c(e_k)c(e_l)c(e_i)\sum_{i\neq
k_1,l_1}2a_{k_1l_1}c(e_i)c(e_{k_1})c(e_{l_1})\\
&=&-4\sum_{i\neq
k,l,k_1,l_1}a_{kl}a_{k_1l_1}c(e_k)c(e_l)c(e_{k_1})c(e_{l_1})\\
&=&-4\left(\sum_{k=k_1,l\neq l_1}+\sum_{k=l_1,l\neq
k_1}+\sum_{l=k_1,k\neq l_1}+\sum_{l=l_1,k\neq
k_1}\right.\\
&&\left.+\sum_{k=k_1,l= l_1}+\sum_{k=l_1,l=
k_1}\right)\left[a_{kl}a_{k_1l_1}c(e_k)c(e_l)c(e_{k_1})c(e_{l_1})\right]\\
&=&-16\sum_{l\neq l_1}a_{kl}a_{kl_1}c(e_l)c(e_{l_1})+16a_{kl}^2.
~~~~~~~~~~~~~~~~~~~~~~~~~~~~~~~~~~~~~~~~~~~~~~~~~~~~(4.22)
\end{eqnarray*}
Similar to (4.22), we have
$$[a_{kl}c(e_k)c(e_l)]^2=\sum_{k\neq l\neq k_1\neq
l_1}a_{kl}a_{k_1l_1}c(e_k)c(e_l)c(e_{k_1})c(e_{l_1})+4\sum_{l\neq
l_1}a_{kl}a_{kl_1}c(e_l)c(e_{l_1})-2a_{kl}^2.\eqno(4.23)$$ Then
$$E=-\frac{s}{4}+2e_k(a_{kl})c(e_l)-\sum_{k\neq l\neq k_1\neq
l_1}a_{kl}a_{k_1l_1}c(e_k)c(e_l)c(e_{k_1})c(e_{l_1})-2a_{kl}^2.\eqno(4.24)$$
Now we can compute ${\rm Tr}(E^2)$.
$${\rm
Tr}[2e_k(a_{kl})c(e_l)2e_{k_1}(a_{k_1l_1})c(e_{l_1})]=-4de_k(a_{kl})e_{k_1}(a_{k_1l})=-4|\delta
\Psi|^2(x_0),\eqno(4.25)$$
where $\delta$ is the adjoint operator of $d$. We have
$${\rm
Tr}[(-\frac{s}{4}-2|\Psi|^2)^2]=4[\frac{s^2}{16}+s|X|^2+4|X|^4],\eqno(4.26)$$
and
$$\sum_{k\neq l\neq k_1\neq
l_1}\sum_{p\neq q\neq r\neq t}{\rm
Tr}\left[a_{kl}a_{k_1l_1}c(e_k)c(e_l)c(e_{k_1})c(e_{l_1})a_{pq}a_{rt}c(e_p)c(e_q)c(e_{r})c(e_{t})\right]$$
$$=\sum_{\{k,l,k_1,l_1\}=\{p,q,r,t\}}{\rm
Tr}\left[a_{kl}a_{k_1l_1}c(e_k)c(e_l)c(e_{k_1})c(e_{l_1})a_{pq}a_{rt}c(e_p)c(e_q)c(e_{r})c(e_{t})\right]$$
$$=8d|\Psi|^4-16da_{kl}a_{k_1l_1}a_{kk_1}a_{ll_1}=32|\Psi|^4-4i_{e_k}i_{e_l}(\Psi)i_{e_{k_1}}i_{e_{l_1}}(\Psi)
i_{e_k}i_{e_{k_1}}(\Psi)i_{e_l}i_{e_{l_1}}(\Psi).\eqno(4.27)$$ By
(4.24)-(4.27), we get
$${\rm Tr}(E^2)=4\left[-|\delta
\Psi|^2+\frac{s^2}{16}+s|\Psi|^2+12|\Psi|^4\right.$$
$$\left.-i_{e_k}i_{e_l}(\Psi)i_{e_{k_1}}i_{e_{l_1}}(\Psi)
i_{e_k}i_{e_{k_1}}(\Psi)i_{e_l}i_{e_{l_1}}(\Psi)\right].\eqno(4.28)$$
\indent In the following, we compute ${\rm
Tr}[\Omega(e_i,e_j)\Omega(e_i,e_j)]$. By (2.14), we have
$$\Omega(e_i,e_j)=R^{S(TM)}(e_i,e_j)-\frac{1}{4}\nabla^{S(TM)}_{e_i}(c(e_j)\Psi+\Psi
c(e_j))$$
$$+\frac{1}{4}\nabla^{S(TM)}_{e_j}(c(e_i)\Psi+\Psi
c(e_i))+\frac{1}{4}\left[c([e_i,e_j])\Psi+\Psi
c([e_i,e_j])\right].\eqno(4.29)$$ Since $\nabla^{S(TM)}$ is a Clifford
connection and $\nabla^{TM}$ has no torsion, we get
$$\Omega(e_i,e_j)=R^{S(TM)}(e_i,e_j)-\frac{1}{4}\nabla^{S(TM)}_{e_i}(\Psi)c(e_j)-\frac{1}{4}c(e_j)\nabla^{S(TM)}_{e_i}(\Psi)
$$
$$+\frac{1}{4}\nabla^{S(TM)}_{e_j}(\Psi)c(e_i)+\frac{1}{4}c(e_i)\nabla^{S(TM)}_{e_j}(\Psi).\eqno(4.30)$$
Since the trace of the odd Clifford elements is zero and
$${\rm Tr}\left[\sum_{i,j=1}^n\left(\frac{1}{4}\sum_{s,t=1}^nR^M_{ijst}c(e_s)c(e_t)\right)^2\right]=-\frac{1}{2}\sum_{i,j,s,t=1}^nR_{ijst}^2,\eqno(4.31)$$
we get
$${\rm
Tr}\sum_{i,j=1}^n[\Omega_{i,j}^2]=-\frac{1}{2}\sum_{i,j,s,t=1}^nR_{ijst}^2+\frac{1}{16}{\rm
Tr}\left[\left.\sum_{i,j=1}^n(-\nabla^{S(TM)}_{e_i}(\Psi)c(e_j)-c(e_j)\nabla^{S(TM)}_{e_i}(\Psi)\right.\right.$$
$$\left.\left.
+\nabla^{S(TM)}_{e_j}(\Psi)c(e_i)+c(e_i)\nabla^{S(TM)}_{e_j}(\Psi)\right)^2\right].\eqno(4.32)$$
By (4.7) (4.28) and (4.32), we obtain\\

\noindent {\bf Proposition 4.2} {\it Let $\Psi$ be a two-form and let $M$ be a $4$-dimensional compact spin manifold without boundary.
 Then}
$$a_4(D_\Psi)=\frac{1}{1440\pi^2}\left\{\Delta(3s+120|\Psi|^2)+\frac{5}{4}s^2-2\sum_{i,l,j,k=1}^nR_{ijik}R_{ljlk}-\frac{7}{4}\sum_{i,j,k,l=1}^nR_{ijkl}^2
+60s|\Psi|^2\right.$$
$$-180|\delta\Psi|^2+2160|\Psi|^4-180\sum_{k,l,k_1,l_1=1}^ni_{e_k}i_{e_l}(\Psi)i_{e_{k_1}}i_{e_{l_1}}(\Psi)
i_{e_k}i_{e_{k_1}}(\Psi)i_{e_l}i_{e_{l_1}}(\Psi)$$
$$
+\frac{15}{8}{\rm
Tr}\left[\sum_{i,j=1}^n\left(-\nabla^{S(TM)}_{e_i}(\Psi)c(e_j)-c(e_j)\nabla^{S(TM)}_{e_i}(\Psi)\right.\right.$$
$$\left.\left.\left.
+\nabla^{S(TM)}_{e_j}(\Psi)c(e_i)+c(e_i)\nabla^{S(TM)}_{e_j}(\Psi)\right)^2\right]\right\}.\eqno(4.33)$$\\

\noindent{\bf Remark.} In fact, Proposition 4.2 holds true, under analogous hypotheses for manifolds of arbitrary dimension and a general two-form $\Psi$ after revising some coefficients.\\

 \noindent {\bf Acknowledgement.} This work
was supported by NCET-13-0721 and NSFC No.11271062.  I would like to thank referees for their careful reading and helpful comments. I also thank Prof. Huitao Feng for his helpful comments.\\

\noindent{\large \bf References}\\

\noindent [Ac] T. Ackermann, A note on the Wodzicki residue, J.
Geom. Phys., 20(1996), 404-406.\\
\noindent[AT]T. Ackermann, J. Tolksdorf, A generalized Lichnerowicz formula, the Wodzicki residue and gravity.
 J. Geom. Phys. 19(1996), no. 2, 143-150.\\
\noindent[BC1]U. Battisti, S. Coriasco, A note on the Einstein-Hilbert action and Dirac operators on $R^n$,
 J. Pseudo-Differ. Oper. Appl. 2(2011),no. 3, 303-315.\\
\noindent[BC2] U. Battisti, S. Coriasco, Wodzicki residue for operators on manifolds with cylindrical ends, Ann. Global Anal. Geom. 40, (2011)
no. 2, 223-249.\\
 \noindent[CC] Chamseddine, A. H.; Connes, A.
Inner fluctuations of the spectral action. J.
Geom. Phys. 57(2006), 1-21.\\
\noindent [Co1] A. Connes, Quantized calculus and applications,
XIth International Congress of Mathematical Physics (Paris,1994),
15-36, Internat Press, Cambridge, MA, 1995.\\
\noindent [Co2] A. Connes, The action functional in
noncommutative
geometry, Comm. Math. Phys., 117(1998),673-683.\\
  \noindent[CM]A. Connes, H. Moscovici, Type III and spectral triples. Traces in Geometry, Number Theory and Quantum
Fields, Aspects of Mathematics E38, Vieweg Verlag 2008, 57-71.\\
\noindent [FGLS] B. V. Fedosov, F. Golse, E. Leichtnam, and E.
Schrohe, The noncommutative residue for manifolds with
boundary, J. Funct.
Anal, 142(1996), 1-31.\\
\noindent[Gi]P. Gilkey, Invariance theory, the heat
equation, and the Atiyah-Singer index theorem, Mathematics Lecture
Series, 11. Publish or Perish, 1984 \\
\noindent [Gu] V.W. Guillemin, A new proof of Weyl's
formula on the asymptotic distribution of eigenvalues, Adv. Math.
55(1985) no.2, 131-160.\\
\noindent [HPS]F. Hanisch, F. Pf\"{a}ffle, C. Stephan, The
spectral action for Dirac operators with skew-symmetric torsion.
Comm. Math. Phys. 300(2010), no. 3, 877-888.\\
\noindent[IL]B. Iochum, C. Levy, Tadpoles and commutative spectral triples. J. Noncommut. Geom. 5 (2011), no. 3, 299-329.\\
\noindent [Ka] D. Kastler, The Dirac operator and gravitation, Commun. Math. Phys, 166(1995), 633-643.\\
\noindent [KW] W. Kalau and M.Walze, Gravity, non-commutative
geometry, and the Wodzicki residue, J. Geom. Phys., 16(1995),327-344.\\
\noindent[Ni]F. Nicola, Trace functionals for a class of pseudo-differential operators
in ${\mathbb{R}}^n$, Math. Phys. Anal. Geom., 6 (2003) no. 1, 89-105.\\
\noindent[Po]R. Ponge, Noncommutative residue for Heisenberg manifolds. Applications in CR and contact geometry,
 J. Funct. Anal. 252(2007), no. 2, 399-463.\\
\noindent [Sc] E. Schrohe, Noncommutative residue, Dixmier's
trace, and heat trace expansions on manifolds with boundary, Contemp. Math. 242(1999), 161-186.\\
\noindent[SZ]A. Sitarz, A. Zajac, Spectral action for scalar perturbations of Dirac operators,
Lett. Math. Phys. 98(2011), no. 3, 333-348.\\
\noindent[Wa1]Y. Wang, Gravity and the noncommutative residue for manifolds with boundary,
 Lett. Math. Phys. 80(2007), 37-56.\\
\noindent[Wa2]Y. Wang, Lower-dimensional volumes and Kastler-Kalau-Walze type theorem for manifolds with boundary,
Commun. Theor. Phys. 54(2010), 38-42.\\
\noindent [Wo] M. Wodzicki, Noncommutative
residue, Chapter I. Fundamentals, Lecture Notes in Math. 1289, 320-399, 1987.\\

 \indent{ School of Mathematics and Statistics,
Northeast Normal University, Changchun Jilin, 130024, China }\\
\indent E-mail: {\it wangy581@nenu.edu.cn}\\

\end{document}